\documentclass{article}
\newcommand{\proof}{\noindent {\bf Proof: }}

\newtheorem{theorem}{Theorem}

\newtheorem{lemma}{Lemma}
\newtheorem{defi}{Definition}

\def\qed{\hfill $\Box$}

\usepackage[]{amssymb}
\usepackage{graphicx}

\begin{document}
\title{Premanifolds}
\author{\'A.G.Horv\'ath}
\date{~}

\maketitle

\begin{abstract}
The tangent hyperplanes of the "manifolds" of this paper equipped a so-called Minkowski product. It is neither symmetric nor bilinear. We give a method to handing such an object as a locally hypersurface of a generalized space-time model and define the main tools of its differential geometry: its fundamental forms, its curvatures and so on. In the case, when the fixed space-time component of the embedding structure is a continuously differentiable semi-inner product space, we get a natural generalization of some important semi-Riemann manifolds as the hyperbolic space, the de Sitter sphere and the light cone of a Minkowski-Lorenz space, respectively.
\end{abstract}



\section{Introduction}

There is no and we will not give a formal definition of an object calling in this paper \emph{ premanifold }. We use this word for a set if it has a manifold-like structure with high freedom in the choosing of the distance function on its tangent hyperplanes. For example we get premanifolds if we investigate the hypersurfaces of a generalized space-time model. The most important types of manifolds as Riemannian, Finslerian or semi-Riemannian can be investigated in this way. The structure of our embedding space was introduced in \cite{gho} and in this paper we shall continue the investigations by the build up the differential geometry of its hypersurfaces. We will give the pre-version of the usual semi-Riemannian or Finslerian spaces, the hyperbolic space, the de Sitter sphere, the light cone and the unit sphere of the rounding semi-inner product space, respectively. In the case, when the space-like component of the generalized space-time model is a continuously differentiable semi-inner product space then we will get back the known and usable geometrical informations on the corresponding hypersurfaces of a pseudo-Euclidean space, e.g. we will show that a prehyperbolic space has constant negative curvature.

\subsection{Terminology}

\begin{description}

\item[concepts without definition:] basis, dimension, direct sum of subspaces, hyperboloid, hyperbolic space and hyperbolic metric, inner (scalar) product, linear and bilinear mapping, real and complex vector spaces, quadratic forms, Riemann, Finsler and semi-Riemann manifolds.
\item[acceleration vector field:] See before Definition 16.
\item[arc-length:] See section 2.2.
\item[Convexity of a hypersurface:] See Definition 10.
\item[Curvature of a curve:] See Definition 14.
\item[de Sitter sphere:] See in paragraph 3.2.
\item[Fundamental forms:] See Definition 11 and 12.
\item[generalized Minkowski space:] See Definition 5.
\item[generalized space-time model:] Finite dimensional, real, generalized Min\-kows\-ki space with one dimensional time-like orthogonal direct components.
\item[geodesic:]See Definition 16.
\item[hypersurface:] The definition in a generalized Minkowski space can be found before Lemma 3.
\item[imaginary unit sphere:] See Definition 8.
\item[i.i.p:] Indefinite inner product (See Definition 3).
\item[Minkowski product:] See Definition 5.
\item[Minkowski-Finsler space:] See Definition 9.
\item[Sectional curvature:] See Definition 15.
\item[s.i.i.p:] Semi-indefinite-inner-product (See Definition 4).
\item[s.i.p:] Semi-inner product (See Definition 1).
\item[Ricci and scalar curvature:] See Definition 16.
\item[tangent vector, tangent hyperplane:] These definitions can be seen before Lemma 3.
\item[velocity vector field:] See before Definition 16.
\end{description}

\subsection{Notation}

\begin{description}

\item[$\mathbb{C}$, $\mathbb{R}$, $\mathbb{R}^{n}$, $S^n$:] The complex line, the real line, the $n$-dimensional
    Euclidean space and the $n$-dimensional unit sphere,
    respectively.
\item[{$\langle \cdot,\cdot\rangle$}:] The notion of scalar product and all its suitable generalizations.
\item[{$[\cdot,\cdot]^-$}:] The notion of s.i.p. corresponding to a generalized Minkowski space.
\item[{$[\cdot,\cdot]^+$}:] The notion of Minkowski product of a generalized Minkowski space.
\item[$f'$:] The derivative of a real-valued function $f$ with domain in $\mathbb{R}$.
\item[$Df$:] The Frechet derivative of a map between two normed spaces.
\item[$f'_e$:] The directional derivative of a real-valued function of a normed space into the direction of $e$.
\item[{$[x,\cdot]'_z(y)$}:] The derivative map of an s.i.p. in its second argument, into the direction of $z$ at the point $(x,y)$. See Definition 3.
\item[$\|\cdot\|'_x(y)$, $\|\cdot\|''_{x,z}(y)$:]The derivative of the norm in the direction of $x$ at the point $y$, and the second derivative of the norm in the directions $x$ and $z$ at the point $y$.
\item[$\Re{\{\cdot\}}$,$\Im{\{\cdot\}}$:] The real and imaginary part of a complex number, respectively.
\item[$T_v$:] The tangent space of a Minkowskian hypersurface at its point $v$.
\item[$\mathcal{S},\mathcal{T},\mathcal{L}$:]The set of space-like, time-like and light-like vectors respectively.
\item[$S$,$T$:]The space-like and time-like orthogonal direct components of a generalized Minkowski space, respectively.
\item[$\{e_1,\ldots,e_k,e_{k+1},\ldots,e_n\}$:] An Auerbach basis of a generalized Minkowski space with $\{e_1,\ldots,e_k\}\subset S$ and $\{e_{k+1},\ldots,e_n\}\subset T$, respectively. All of the $e_i'$ orthogonal to the another ones with respect to the Minkowski product.
\item[$G$,$G^+$:] The unit sphere of a generalized space-time model and its upper sheet, respectively.
\item[$H$,$H^+$:] The sphere of radius $i$ and its upper sheet, respectively.
\item[$K$, $K^+$:] The unit sphere of the embedding semi-inner product space and its upper sheet, respectively.
\item[$L$,$L^+$:] The light cone of a generalized space-time model and its upper sheet, respectively.
\item[$g$:] The function $g(s)=s+\mathfrak{g}(s)e_n$ with $\mathfrak{g}(s)=\sqrt{-1+[s,s]}$ defines the points of $G:=\{s+g(s) | s\in S$\}.
\item[$h$:] The function $h(s)=s+\mathfrak{h}(s)e_n$ with $\mathfrak{h}(s)=\sqrt{1+[s,s]}$ defines the points of $H^+:=\{s+h(s) | s\in S$\}.
\item[$k$:] The function $k(s)=s+\mathfrak{k}(s)e_n$ with $\mathfrak{k}(s)=\sqrt{1-[s,s]}$ defines the points of $K^+:=\{s+k(s) | s\in S$\}.
\item[$l$:] The function $l(s)=s+\mathfrak{l}(s)e_n$ with $\mathfrak{l}(s)=\sqrt{[s,s]}$ defines the points of $L^+:=\{s+l(s) | s\in S$\}.
\end{description}

\subsection{History with completion of the preliminaries}

A generalization of the inner product and the inner product spaces (briefly i.p spaces) was raised by G. Lumer in \cite{lumer}.
\begin{defi}[\cite{lumer}]
The \emph{ semi-inner-product (s.i.p)} on a complex vector space $V$ is a complex function $[x,y]:V\times V\longrightarrow \mathbb{C}$ with the following properties:
\begin{description}

\item[s1]: $[x+y,z]=[x,z]+[y,z]$,
\item[s2]: $[\lambda x,y]=\lambda[x,y]$ \mbox{ for every } $\lambda \in \mathbb{C}$,
\item[s3]: $[x,x]>0$ \mbox{ when } $x\not =0$,
\item[s4]: $|[x,y]|^2\leq [x,x][y,y]$.
\end{description}
A vector space $V$ with a s.i.p. is an \emph{ s.i.p. space}.
\end{defi}

G. Lumer proved that an s.i.p space is a normed vector space with norm $\|x\|=\sqrt{[x,x]}$ and, on the other hand, that every normed vector space can be represented as an s.i.p. space. In \cite{giles} J. R. Giles showed that the following homogeneity property holds:

\begin{description}

\item[s5]: $[x,\lambda y]=\bar{\lambda}[x,y]$ for all complex $\lambda $.

\end{description}
This can be imposed, and all normed vector spaces can be represented as s.i.p. spaces with this property. Giles also introduced the concept of {\bf continuous s.i.p. space} as an s.i.p. space having the additional property

\begin{description}
\item[s6]: For any unit vectors $x,y \in S$, $\Re\{[y,x+\lambda y]\}\rightarrow\Re\{[y,x]\}$ for all real $\lambda\rightarrow 0$.
\end{description}
The space is uniformly continuous if the above limit is reached uniformly for all points $x,y$ of the unit sphere $S$.
A characterization of the continuous s.i.p. space is based on the differentiability property of the space.

Giles proved in \cite{giles} that
\begin{theorem}[\cite{giles}]
An s.i.p. space is a continuous (uniformly continuous) s.i.p. space if and only if the norm is G\^{a}teaux (uniformly Fr\`{e}chet) differentiable.
\end{theorem}

In \cite{gho} \'A.G.Horv\'ath defined the differentiable s.i.p. as follows:

\begin{defi}
A \emph{ differentiable s.i.p. space} is an continuous s.i.p. space where the s.i.p. has the additional property

\noindent {\bf s6':} For every three vectors x,y,z and real $\lambda $
$$
[x,\cdot]'_z(y):=\lim\limits_{\lambda \rightarrow 0}\frac{\Re\{[x,y+\lambda z]\}-\Re\{[x,y]\}}{\lambda }
$$
does exist. We say that the s.i.p. space is \emph{ continuously differentiable}, if the above limit, as a function of $y$, is continuous.
\end{defi}

First we note that the equality $\Im\{[x,y]\}=\Re\{[-ix,y]\}$ together with the above property guarantees the existence and continuity of the complex limit:
$$
\lim\limits_{\lambda \rightarrow 0}\frac{[x,y+\lambda z]-[x,y]}{\lambda }.
$$

The following theorem was mentioned without proof in \cite{gho}:

\begin{theorem}[\cite{gho}]
An s.i.p. space is a (continuously) differentiable s.i.p. space if and only if the norm is two times (continuously) G\^{a}teaux differentiable. The connection between the derivatives is
$$
\|y\|(\|\cdot\|''_{x,z}(y))= [x,\cdot]'_z(y)-\frac{\Re[x,y]\Re[z,y]}{\|y\|^2}.
$$
\end{theorem}

Since the present paper often use this statement, we give a proof for it. We need the following useful lemma going back, with different notation, to McShane \cite{mcshane} or Lumer \cite{lumer2}.

\begin{lemma}[\cite{lumer2}]
If E is any s.i.p. space with $x,y\in E$, then
$$
\|y\| (\|\cdot\|'_x(y))^-\leq \Re\{[x,y]\}\leq \|y\| (\|\cdot\|'_x(y))^+
$$
holds, where $(\|\cdot\|'_x(y))^-$ and $(\|\cdot\|'_x(y))^+$ denotes the left hand and right hand derivatives with respect to the real variable $\lambda $. In particular, if the norm is differentiable, then
$$
[x,y]= \|y\| \{(\|\cdot\|'_x(y))+\|\cdot\|'_{-ix}(y)\}.
$$
\end{lemma}

Now we prove Theorem 2.

\proof[of Theorem 2]
To determine the derivative of the s.i.p., assume that the norm is twice differentiable. Then, by Lemma 1 above, we have
$$
\frac{\Re\{[x,y+\lambda z]\}-\Re\{[x,y]\}}{\lambda }=\frac{\|y+\lambda z\|(\|\cdot\|'_x(y+\lambda z))-\|y\| (\|\cdot\|'_x(y))}{\lambda }=
$$
$$
=\frac{\|y\|\|y+\lambda z\|(\|\cdot\|'_x(y+\lambda z))-\|y\|^2 (\|\cdot\|'_x(y))}{\lambda \|y\|}\geq
$$
$$
\geq \frac{|[y+\lambda z,y]|(\|\cdot\|'_x(y+\lambda z))-\|y\|^2 (\|\cdot\|'_x(y))}{\lambda \|y\|},
$$
 where we have assumed that the sign of $\frac{\|\cdot\|'_x(y+\lambda z)}{\lambda}$ is positive. Since the derivative of the norm is continuous, this follows from the assumption that $\frac{\|\cdot\|'_x(y)}{\lambda}$ is positive. Considering the latter condition, we get
$$
\frac{\Re\{[x,y+\lambda z]\}-\Re\{[x,y]\}}{\lambda }\geq
$$
$$
\geq \|y\|^2 \frac{\|\cdot\|'_x(y+\lambda z)-(\|\cdot\|'_x(y))}{\lambda \|y\|}+\frac{\Re[z,y]}{\|y\|}\|\cdot\|'_x(y+\lambda z).
$$
On the other hand,
$$
\frac{\|y+\lambda z\|(\|\cdot\|'_x(y+\lambda z))-\|y\| (\|\cdot\|'_x(y))}{\lambda }\leq
$$
$$
\leq \frac{\|y+\lambda z\|^2(\|\cdot\|'_x(y+\lambda z))-|[y,y+\lambda z]| (\|\cdot\|'_x(y))}{\lambda \|y+\lambda z\|}=
$$
$$
=\frac{\|y+\lambda z\|^2(\|\cdot\|'_x(y+\lambda z))- (\|\cdot\|'_x(y))}{\lambda \|y+\lambda z\|}+\lambda \Re{[z,y+\lambda z]}\frac{(\|\cdot\|'_x(y))}{\lambda \|y+\lambda z\|}.
$$
Analogously, if $\frac{\|\cdot\|'_x(y)}{\lambda}$ is negative, then both of the above inequalities are reversed, and we get that the limit
$$
\lim\limits_{\lambda \mapsto 0}\frac{\Re\{[x,y+\lambda z]\}-\Re\{[x,y]\}}{\lambda }
$$
exists, and equals to
$$
\|y\|(\|\cdot\|''_{x,z}(y))+\frac{\Re[x,y]\Re[z,y]}{\|y\|^2}.
$$
Here we note that also in the case $\frac{\|\cdot\|'_x(y)}{\lambda}=0$ there exists a neighborhood in which the sign of the function $\frac{\|\cdot\|'_x(y+\lambda z)}{\lambda}$ is constant. Thus we, need not investigate this case by itself.
Conversely, consider the fraction
$$
\|y\|\frac{\|\cdot\|'_x(y+\lambda z)-(\|\cdot\|'_x(y))}{\lambda}.
$$
We assume now that the s.i.p. is differentiable, implying that it is continuous, too. The norm is differentiable by the theorem of Giles. Using again Lemma 1 and  assuming that $\frac{\Re[x,y]}{\lambda}> 0$, we have
$$
\|y\|\frac{\|\cdot\|'_x(y+\lambda z)-(\|\cdot\|'_x(y))}{\lambda}=\frac{\Re[x,y+\lambda z]\|y\|-\Re[x,y]\|y+\lambda z\|}{\lambda \|y+\lambda z\|}=
$$
$$
=\frac{\Re[x,y+\lambda z]\|y\|^2-\Re[x,y]\|y+\lambda z\|\|y\|}{\lambda \|y\|\|y+\lambda z\|}\leq
$$
$$
\frac{\Re[x,y+\lambda z]\|y\|^2-\Re[x,y]|[y+\lambda z,y]|}{\lambda \|y\|\|y+\lambda z\|}=
$$
$$
=\frac{\Re\{[x,y+\lambda z]\}-\Re\{[x,y]\}}{\lambda }\frac{\|y\|}{\|y+\lambda z\|}-\frac{\Re[x,y]\Re[z,y]}{\|y\|\|y+\lambda z\|}.
$$
On the other hand, using the continuity of the s.i.p. and our assumption $\frac{\Re[x,y]}{\lambda}> 0$ similarly as above, we also get an inequality:
$$
\|y\|\frac{\|\cdot\|'_x(y+\lambda z)-(\|\cdot\|'_x(y))}{\lambda}\geq
$$
$$
\frac{\Re\{[x,y+\lambda z]\}-\Re\{[x,y]\}}{\lambda} -\frac{\Re[x,y+\lambda z]\Re[z,y+\lambda z]}{\|y+\lambda z\|^2}.
$$
If we reverse the assumption of signs, then the direction of the inequalities will also change. Again a limit argument shows that the first differential function is differentiable, and the connection between the two derivatives is
$$
\|y\|(\|\cdot\|''_{x,z}(y))= [x,\cdot]'_z(y)-\frac{\Re[x,y]\Re[z,y]}{\|y\|^2}.
$$
\qed

From geometric point of view we know that if $C$ is a $0$-symmetric, bounded,
convex body in the Euclidean $n$-space $\mathbb{R}^n$ (with fixed
origin O), then it defines a norm whose unit ball is $C$ itself (see
\cite{l-g}). Such a space is called (Minkowski or) normed linear space. Basic
results on such spaces are collected in the surveys
\cite{martini-swanepoel 1}, \cite{martini-swanepoel 2}, and
\cite{martini}. In fact, the norm is a continuous function which is
considered (in geometric terminology, as in \cite{l-g}) as a gauge
function. Combining  this with the result of Lumer and Giles we get
that a normed linear space can be represented as an s.i.p space.
The metric of such a space (called Minkowski metric), i.e., the distance of
two points induced by this norm, is invariant with respect to translations.

Another concept of  Minkowski space was also raised by H. Minkowski and used in Theoretical Physics and Differential Geometry, based on the
concept of indefinite inner product. (See, e.g., \cite{gohberg}.)

\begin{defi}[\cite{gohberg}]
The \emph{indefinite inner product (i.i.p.)} on a complex vector space $V$
is a complex function
$
[x,y]:V\times V\longrightarrow \mathbb{C}
$
with the following properties:
\begin{description}

\item[i1]: $[x+y,z]=[x,z]+[y,z]$,
\item[i2]: $[\lambda x,y]=\lambda[x,y]$ \mbox{ for every } $\lambda \in \mathbb{C}$,
\item[i3]: $[x,y]=\overline{[y,x]}$ \mbox{ for every } $x,y\in V$,
\item[i4]: $[x,y]=0$ \mbox{ for every } $y\in V$ then $x=0$.
\end{description}
A vector space $V$ with an i.i.p. is called an \emph{ indefinite inner product space}.
\end{defi}

The standard  mathematical model of space-time is a four
dimensional i.i.p. space with signature $(+,+,+,-)$, also
called Minkowski space in the literature. Thus we have a well known homonymism with the notion of Minkowski space!

In \cite{gho} the concepts of s.i.p. and i.i.p. was combined in the following one:

\begin{defi}[\cite{gho}]
The \emph{ semi-indefinite inner product (s.i.i.p.)} on a complex vector space $V$ is a complex function $[x,y]:V\times V\longrightarrow \mathbb{C}$ with the following properties:
\begin{description}

\item[1] $[x+y,z]=[x,z]+[y,z]$ (additivity in the first argument),
\item[2] $[\lambda x,y]=\lambda[x,y]$ \mbox{ for every } $\lambda \in \mathbb{C}$ (homogeneity in the first argument),
\item[3] $[x,\lambda y]=\overline{\lambda}[x,y]$ \mbox{ for every } $\lambda \in \mathbb{C}$ (homogeneity in the second argument),
\item[4] $[x,x]\in \mathbb{R}$ \mbox{ for every } $x\in V$ (the corresponding quadratic form is real-valued),
\item[5] if either $[x,y]=0$ \mbox{ for every } $y\in V$ or $[y,x]=0$ for all $y\in V$, then $x=0$ (nondegeneracy),
\item[6] $|[x,y]|^2\leq [x,x][y,y]$ holds on non-positive and non-negative subspaces of V, respectively (the Cauchy-Schwartz inequality is valid on positive and negative subspaces, respectively).
\end{description}
A vector space $V$ with an s.i.i.p. is called an \emph{ s.i.i.p. space}.
\end{defi}

It was conclude that an s.i.i.p. space is a
homogeneous s.i.p. space if and only if the property {\bf s3} holds,
too. An s.i.i.p. space is an i.i.p. space if and only if the
s.i.i.p. is an antisymmetric product. In this latter case
$[x,x]=\overline{[x,x]}$ implies {\bf 4}, and the function is
also Hermitian linear in its second argument. In fact, we have:
$[x,\lambda y+\mu z]=\overline{[\lambda y+\mu
z,x]}=\overline{\lambda} \overline{[y,x]}+\overline{\mu} \overline{[
z,x]}=\overline{\lambda}[x,y]+\overline{\mu}[x,z]$. It is clear that
both of the classical "Minkowski spaces" can be represented either by an
s.i.p or by an i.i.p., so automatically they can also be represented as an
s.i.i.p. space.

The following fundamental lemma was proved in \cite{gho}:

\begin{lemma}[\cite{gho}]
Let  $(S,[\cdot,\cdot]_S)$ and $(T,-[\cdot,\cdot]_T)$ be two s.i.p.
spaces. Then the function
$[\cdot,\cdot]^-:(S+T)\times(S+T)\longrightarrow \mathbb{C}$ defined
by
$$
[s_1+t_1,s_2+t_2]^-:=[s_1,s_2]-[t_1,t_2]
$$
is an s.i.p. on the vector space $S+T$.
\end{lemma}

It is possible that the s.i.i.p. space $V$ is a direct sum of its two subspaces where one of them is positive and the other one is negative. Then there are two more structures on $V$, an s.i.p. structure (by Lemma 2) and a natural third one, which was called by Minkowskian  structure.
\begin{defi}[\cite{gho}]
Let $(V,[\cdot,\cdot])$ be an s.i.i.p. space. Let $S,T\leq V$ be positive and negative subspaces, where $T$ is a direct complement of $S$ with
respect to $V$. Define a product on $V$ by the equality $[u,v]^+=[s_1+t_1,s_2+t_2]^+=[s_1,s_2]+[t_1,t_2]$, where $s_i\in S$ and $t_i\in T$,
respectively.  Then we say that the pair $(V,[\cdot,\cdot]^+)$ is a \emph{ generalized Minkowski space with Minkowski product $[\cdot,\cdot]^+$}. We also say
that $V$ is a \emph{ real generalized Minkowski space} if it is a real vector space and the s.i.i.p. is a real valued function.
\end{defi}

The Minkowski product defined by the above equality satisfies properties {\bf 1}-{\bf 5} of the s.i.i.p.. But in general, property {\bf 6} does not hold. (See an example in \cite{gho}.)

By Lemma 2 the function $\sqrt{[v,v]^-}$ is a norm function on $V$ which can give an embedding space for a generalized Minkowski
space. This situation is analogous to the situation when a pseudo-Euclidean space is obtained from a Euclidean space by the action of an i.i.p.

It is easy to see that by the methods of \cite{gho}, starting with arbitrary two normed spaces $S$ and $T$, one can mix a generalized Minkowski space. Of course its smoothness property is basically determined by the analogous properties of $S$ and $T$.

If now we consider the theory of s.i.p in the sense of Lumer-Giles, we have a natural concept of orthogonality. For the unified terminology we change the original notation of Giles and use instead

\begin{defi}[\cite{giles}]
The vector \emph{ $x$ is orthogonal to the vector $y$} if $[x,y]=0$.
\end{defi}
Since s.i.p. is neither antisymmetric in the complex case nor symmetric in the real one, this definition of orthogonality is not symmetric in general.

Giles proved that in a continuous s.i.p. space $x$ is orthogonal to $y$ in the sense of the s.i.p. if and only if $x$ is orthogonal to $y$ in the sense of Birkhoff-James. (See e.g. \cite{alonso1} and \cite{alonso2}.) We note that the s.i.p. orthogonality implies the Birkhoff-James orthogonality in every normed space. Lumer pointed out that a normed linear space can be transformed into an s.i.p. space in a unique way if and only if its unit sphere is smooth (i.e., there is a unique supporting hyperplane at each point of the unit sphere). In this case the corresponding (unique) s.i.p. has the homogeneity property {\bf [s5]}.

Let $(V,[\cdot,\cdot])$ be an s.i.i.p. space, where $V$ is a complex (real) vector space. It was defined the orthogonality of such a space by a definition analogous to the definition of the orthogonality of an i.i.p. or s.i.p. space.
\begin{defi}[\cite{gho}]
The vector \emph{ $v$ is orthogonal to the vector $u$} if $[v,u]=0$. If $U$ is a subspace of $V$, define the orthogonal companion of $U$ in $V$ by
$$
U^\bot =\{v\in V | [v,u]=0 \mbox{ for all } u\in U\}.
$$
\end{defi}

We note that, as in the i.i.p. case, the orthogonal companion is always a subspace of $V$. It was proved the following theorem:

\begin{theorem}[\cite{gho}] Let $V$ be an $n$-dimensional s.i.i.p. space. Then the orthogonal companion of a non-neutral vector $u$ is a subspace having a direct complement of the linear hull of $u$ in $V$. The orthogonal companion of a neutral vector $v$ is a degenerate subspace of dimension $n-1$ containing $v$.
\end{theorem}

Observe that this proof does not use the property {\bf 6} of the
s.i.i.p.. So this statement is true for any concepts of product
satisfying properties {\bf 1}-{\bf 5}. As we saw, the Minkowski product is also such a product.

We also note that in a generalized Minkowski space, the positive and negative components $S$ and $T$ are Pythagorean orthogonal to each other. In fact, for every pair of vectors $s\in S$ and $t\in T$, by definition we have
$[s-t,s-t]^+=[s,s]+[-t,-t]=[s,s]^++[t,t]^+$.

Let $V$ be a generalized Minkowski space. Then we call a vector {\bf space-like, light-like, or time-like} if its scalar square is positive, zero, or
negative, respectively. Let $\mathcal{S}, \mathcal{L}$ and $\mathcal{T}$ denote the sets of the space-like, light-like, and time-like vectors,
respectively.

In a finite dimensional, real generalized Minkowski space with $\dim T=1$ it can geometrically characterize these sets of vectors. Such a space
is called in \cite{gho} a {\bf generalized space-time model}. In this case $\mathcal{T}$ is a union of its two parts, namely
$$
 \mathcal{T}=\mathcal{T}^+\cup \mathcal{T}^-,
$$
where
$$
\mathcal{T}^+=\{s+t\in \mathcal{T} | \mbox{ where } t=\lambda e_n \mbox{ for } \lambda \geq 0\} \mbox{ and }
$$
$$\mathcal{T}^-=\{s+t\in \mathcal{T} |
\mbox{ where } t=\lambda e_n \mbox{ for } \lambda \leq 0\}.
$$

It has special interest, the imaginary unit sphere of a finite dimensional, real, generalized space-time model. (See Def.8 in \cite{gho}.) It was given a metric on it, and thus got a structure similar to the hyperboloid model of the hyperbolic space embedded in a
space-time model. In the case when the space $S$ is an Euclidean space this hypersurface is a model of the $n$-dimensional hyperbolic space thus it is such-like generalization of it.

It was proved in \cite{gho} the following:

\begin{theorem}[\cite{gho}]
Let $V$ be a generalized space-time model. Then $\mathcal{T}$ is an open double cone with boundary $\mathcal{L}$, and the positive part $\mathcal{T}^+$ (resp. negative part $\mathcal{T}^-$) of $\mathcal{T}$ is convex.
\end{theorem}

We note that if $\dim T> 1$ or the space is complex, then the set of time-like vectors cannot be divided into two convex components. So we have to consider that our space is a generalized space-time model.
\begin{defi}[\cite{gho}]
The set
$$
H:=\{ v\in V | [v,v]^+=-1\},
$$
is called the \emph{ imaginary unit sphere } of the generalized space-time model.
\end{defi}

With respect to the embedding real normed linear space
$(V,[\cdot,\cdot]^-)$ (see Lemma 2) $H$ is, as we saw, a generalized two sheets
hyperboloid corresponding to the two pieces of $\mathcal{T}$,
respectively. Usually we deal only with one sheet of the hyperboloid,
or identify the two sheets projectively. In this case the space-time
component $s\in S$ of $v$ determines uniquely the time-like  one, namely
$t\in T$. Let $v\in H$ be arbitrary. Let $T_v$ denote the set
$v+v^\bot$, where $ v^\bot$ is the orthogonal complement of
$v$ with respect to the s.i.i.p., thus a subspace.

It was also proved that the set $T_v$ corresponding to the point $v=s+t\in H$ is a positive, (n-1)-dimensional affine subspace of the generalized Minkowski space
$(V,[\cdot,\cdot]^+)$.

Each of the affine spaces $T_v$ of $H$ can be considered as a semi-metric space, where the semi-metric arises from the Minkowski product restricted
to this positive subspace of $V$. We recall that the Minkowski product does not satisfy the Cauchy-Schwartz inequality. Thus the corresponding distance function does not satisfy the triangle inequality. Such a distance function is called in the literature semi-metric (see \cite{tamassy}). Thus, if the set $H$ is sufficiently smooth, then a metric can be adopted for it, which arises from the restriction of the Minkowski product to the tangent spaces of $H$. Let us discuss this more precisely.

The directional derivatives of a function $\mathfrak{f}:S\longmapsto \mathbb{R}$ with respect to a unit vector $e$ of $S$ can be defined in the usual way, by the existence of the limits for real $\lambda $:
$$
\mathfrak{f}'_{e}(s)=\lim\limits_{\lambda \mapsto 0}\frac{\mathfrak{f}(s+\lambda e)-\mathfrak{f}(s)}{\lambda}.
$$
Let now the generalized Minkowski space be a generalized space-time model, and consider a mapping $f$ on $S$ to $\mathbb{R}$. Denote by $e_n$ a basis vector of $T$ with length $i$ as in the definition of $\mathcal{T}^+$ before Theorem 4. The set of points
$$
F:=\{(s+\mathfrak{f}(s)e_n)\in V \mbox{ for } s\in S\}
$$
is a so-called {\bf hypersurface} of this space. Tangent vectors of a hypersurface $F$ in a point $p$ are the vectors associated to the directional derivatives of the coordinate functions in the usual way. So
$u$ is a {\bf tangent vector} of the hypersurface $F$ in its point $v=(s+\mathfrak{f}(s)e_n)$, if it is of the form
$$
u=\alpha (e+\mathfrak{f}'_{e}(s)e_n) \mbox{ for real } \alpha \mbox{ and unit vector } e\in S.
$$
The linear hull of the tangent vectors translated into the point $s$ is the tangent space of $F$ in $s$. If the tangent space has dimension $n-1$, we call it {\bf tangent hyperplane}.

We now reformulate Lemma 3 of \cite{gho}:
\begin{lemma}[See also in \cite{gho} as Lemma 3]
Let $S$ be a continuous (comp\-lex) s.i.p. space. (So the property {\bf s6} holds.) Then the directional derivatives of the real valued function
$$
\mathfrak{h}:s\longmapsto \sqrt{1+[s,s]}
$$
are
$$
\mathfrak{h}'_{e}(s)=\frac{\Re{[e,s]}}{\sqrt{1+[s,s]}}.
$$
\end{lemma}

The following theorem is a consequence of this result.

\begin{theorem}
Let assume that the s.i.p. $[\cdot,\cdot]$ of $S$ is differentiable. (So the property {\bf s6'} holds.) Then for every two vectors $x$ and $z$ in $S$ we have:
$$
[x,\cdot]'_z(x)=2\Re[z,x]-[z,x],
$$
and
$$
\|\cdot\|''_{x,z}(x)=\frac{\Re[z,x]-[z,x]}{\|x\|}.
$$
If we also assume that the s.i.p. is continuously differentiable (so the norm is a $C^2$ function), then we also have
$$
[x,\cdot]'_x(y)=[x,x],
$$
and thus
$$
\|\cdot\|''_{x,x}(y)=\|x\|^2-\frac{\Re[x,y]^2}{\|y\|^2}.
$$
\end{theorem}

\proof
Since
$$
\frac{1}{\lambda}\left([x+\lambda z,x+\lambda z]-[x,x]\right)=
\frac{1}{\lambda}\left([x,x+\lambda z]-[x,x]\right)+\frac{1}{\lambda}[\lambda z,x+\lambda z],
$$
if $\lambda $ tends to zero then the right hand side tends to
$$
[x,\cdot]'_z(x)+[z,x].
$$
The left hand side is equal to
$$
\frac{\left(\sqrt{1+[x+\lambda z,x+\lambda z]}-\sqrt{1+[x,x]}\right)\left(\sqrt{1+[x+\lambda z,x+\lambda z]}+\sqrt{1+[x,x]}\right)}{\lambda}
$$
thus by Lemma 3 it tends to
$$
\frac{\Re{[z,x]}}{\sqrt{1+[x,x]}}2\sqrt{1+[x,x]}.
$$
This implies the first equality
$$
[x,\cdot]'_z(x)=2\Re[z,x]-[z,x].
$$
Using Theorem 2 in \cite{gho} we also get that
$$
\|x\|(\|\cdot\|''_{x,z}(x))= [x,\cdot]'_z(x)-\frac{\Re[x,x]\Re[z,x]}{\|x\|^2},
$$
proving the second statement, too.

If we assume that the norm is a $C^2$ function of its argument then the first derivative of the second argument of the product is a continuous function of its arguments. So the function $A(y):S\longrightarrow \mathbb{R}$ defined by the formula
$$
A(y)=[x,\cdot]'_x(y)= \lim\limits_{\lambda \mapsto 0}\frac{1}{\lambda}\left([x,y+\lambda x]-[x,y]\right)
$$
continuous in $y=0$. On the other hand for non-zero $t\in \mathbb{R}$ we use the notation $t\lambda'=\lambda $ and we get that
$$
A(ty)=\lim\limits_{\lambda \mapsto 0}\frac{1}{\lambda}\left([x,ty+\lambda x]-[x,y]\right)=\lim\limits_{\lambda' \mapsto 0}\frac{t}{t\lambda '}\left([x,y+\lambda' x]-[x,y]\right)=A(y).
$$
From this we can see immediately that
$$
[x,\cdot]'_x(y)=A(y)=A(0)=[x,x]
$$
holds for every $y$.
Applying again the formula connected the derivative of the product and the norm we get the last statement of the theorem, too.
\qed

Applying Lemma 3 to $H^+$ it was given a connection between the differentiability properties and the orthogonality one.
The tangent vectors of the hypersurface $H^+$ in its point
$$
v=s+\sqrt{1+[s,s]}e_n
$$
form the orthogonal complement $v^{\bot}$ of $v$ with respect to the Minkowski product.

It was defined in \cite{gho} a Finsler space type structure for a hypersurface of a generalized space-time model.

\begin{defi}[\cite{gho}]

Let $F$ be a hypersurface of a generalized space-time model for which the following properties hold:
\begin{description}
\item i, In every point $v$ of $F$, there is a (unique) tangent hyperplane $T_v$ for which the restriction of the Minkowski product
$[\cdot,\cdot]^+_v$ is positive, and

\item ii, the function $ds^2_v:=[\cdot,\cdot]^+_v:F\times T_v\times T_v\longrightarrow \mathbb{R^+}$
$$
ds^2_v:(v,u_1,u_2)\longmapsto [u_1,u_2]^+_v
$$
varies differentiably with the vectors $v\in F$ and $u_1,u_2\in T_v$.
\end{description}
Then we say that the pair \emph{ $(F, ds^2)$ is a Minkowski-Finsler space} with semi-metric $ds^2$ embedding into the generalized space-time model $V$.
\end{defi}

Naturally "varies differentiably with the vectors $v,u_1,u_2$" means that for every $v\in T$ and pairs of vectors $u_1,u_2\in T_v$ the function $[u_1,u_2]_v$ is a differentiable function on $F$.
One of the important results on the imaginary unit sphere was

\begin{theorem}[\cite{gho}]
Let $V$ be a generalized space-time model. Let $S$ be a continuously differentiable s.i.p. space, then $(H^+,ds^2)$ is a Minkowski-Finsler space.
\end{theorem}

In present paper we will prefer the name "pre-hyperbolic space" for this structure.

\begin{center}
{\large \bf Acknowledgment}
\end{center}

The author wish to thank for {\bf G.Moussong} who suggested the investigation of $H^+$ by the tools of differential geometry and {\bf B.Csik\'os} who also gave helpful hints.

\section{ Hypersurfaces as premanifolds}

\subsection{Convexity, fundamental forms}

Let $S$ be a continuously differentiable s.i.p. space,  $V$ be a generalized space-time model and $F$ a hypersurface. We shall say that $F$ is a {\bf space-like hypersurface} if the Minkowski product is positive on its all tangent hyperplanes. The objects of our examination are the convexity, the fundamental forms, the concepts of curvature, the arc-length and the geodesics. In this section we in a generalized space-time model define these that would be a generalizations of the known concepts. In a pseudo-Euclidean or semi-Riemann space it can be found in the notes \cite{moussong} and the book \cite{dubrovin}.

\begin{defi}[\cite{moussong}]
We say that a hypersurface is \emph{convex} if it lies on one side of its each tangent hyperplanes. It is \emph{strictly convex} if it is convex and its tangent hyperplanes contain precisely one points of the hypersurface, respectively.
\end{defi}

In an Euclidean space the {\bf first fundamental form} is a positive definite quad\-ra\-tic form induced by the inner product of the tangent space.

In our generalized space-time model the first fundamental form is  giving by the scalar square of the tangent vectors with respect to the Minkowski product restricted to the tangent hyperplane.
If we have a map $f:S\longrightarrow V$ then it can be decomposed to a sum of its space-like and time-like components. We have
$$
f=f_S+f_T
$$
where  $f_S:S\longrightarrow S$ and $f_T:S\longrightarrow T$, respectively.
With respect to the embedding normed space we can compute its Frechet derivative by the low
$$
Df=\left[\begin{array}{c}
Df_S \\
Df_T
\end{array}\right]
$$
implying that
$$
Df(s)=Df_S(s)+Df_T(s).
$$

Introduce the notation
$$
{[f_1(c(t)),\cdot]^+}'_{D(f_2\circ c)(t)}(f_2(c(t))):=
$$
$$
:=\left([(f_1)_S(c(t)),\cdot]'_{D((f_2)_S\circ c)(t)}((f_2)_S(c(t))) - (f_1)_T(c(t))((f_2)_T\circ c)'(t)\right).
$$

We need the following technical lemma:

\begin{lemma}
If $f_1, f_2: S\longrightarrow V$ are two $C^2$ maps and $c:\mathbb{R}\longrightarrow S$ is an arbitrary $C^2$ curve then
$$
([(f_1\circ c)(t)),(f_2\circ c)(t))]^+)'=
$$
$$
=[D(f_1\circ c)(t),(f_2\circ c)(t))]^++{[(f_1\circ c)(t)),\cdot]^+}'_{D(f_2\circ c)(t)}((f_2\circ c)(t)).
$$
\end{lemma}

\proof
By definition
$$
([f_1\circ c,f_2\circ c)]^+)'|_t:=\lim\limits_{\lambda \rightarrow 0}\frac{1}{\lambda}\left([f_1(c(t+\lambda)),f_2(c(t+\lambda ))]^+-[f_1(c(t)),f_2(c(t))]^+\right)
$$
$$
=\lim\limits_{\lambda \rightarrow 0}\frac{1}{\lambda}\left([(f_1)_S(c(t+\lambda)),(f_2)_S(c(t+\lambda ))]-[(f_1)_S(c(t)),(f_2)_S(c(t))]\right)+
$$
$$
+\lim\limits_{\lambda \rightarrow 0}\frac{1}{\lambda}\left([(f_1)_T(c(t+\lambda)),(f_2)_T(c(t+\lambda))] -[(f_1)_T(c(t)),(f_2)_T(c(t))]\right).
$$
The first part is
$$
\lim\limits_{\lambda \rightarrow 0}\frac{1}{\lambda}\left([(f_1)_S(c(t+\lambda))-(f_1)_S(c(t)),(f_2)_S(c(t+\lambda ))]+\right.
$$
$$
\left.+[(f_1)_S(c(t)),(f_2)_S(c(t+\lambda))]-[(f_1)_S(c(t)),(f_2)_S(c(t))]\right)=
$$
$$
=[D((f_1)_S\circ c)|_{t},(f_2)_S(c(t))]+ [(f_1)_S(c(t)),\cdot]'_{D((f_2)_S\circ c)(t)}((f_2)_S(c(t))).
$$
To prove this take a coordinate system $\{e_1,\cdots ,e_{n-1}\}$ in $S$ and consider the coordinate-wise representation
$$
(f_2)_S\circ c=\sum\limits_{i=1}^{n-1}((f_2)_S\circ c)_ie_i
$$
of $(f_2)_S\circ c$. Using Taylor's theorem for the coordinate functions we have that there are real parameters $t_i \in (t,t+\lambda)$, for which
$$
((f_2)_S\circ c)(t+\lambda)=((f_2)_S\circ c)(t)+\lambda D((f_2)_S\circ c)(t)+\frac{1}{2}\lambda ^2\sum\limits_{i=1}^{n-1}((f_2)_S\circ c)''_i(t_i)e_i.
$$
Thus we can get
$$
[(f_1)_S(c(t)),(f_2)_S(c(t+\lambda))]-[(f_1)_S(c(t)),(f_2)_S(c(t))]=
$$
$$
=[(f_1)_S(c(t)),(f_2)_S(c(t))+D((f_2)_S\circ c)(t)\lambda +
$$
$$
+\frac{1}{2}\lambda ^2\sum\limits_{i=1}^{n-1}((f_2)_S\circ c)''_i(t_i)e_i]-[(f_1)_S(c(t)),(f_2)_S(c(t))]=
$$
$$
\left([(f_1)_S(c(t)),(f_2)_S(c(t))+D((f_2)_S\circ c)(t)\lambda ]-[(f_1)_S(c(t)),(f_2)_S(c(t))]\right)+
$$
$$
+[(f_1)_S(c(t)),(f_2)_S(c(t))+D((f_2)_S\circ c)(t)\lambda +\frac{1}{2}\lambda ^2\sum\limits_{i=1}^{n-1}((f_2)_S\circ c)''_i(t_i)e_i ]-
$$
$$
-[(f_1)_S(c(t)),(f_2)_S(c(t))+D((f_2)_S\circ c)(t)\lambda ].
$$
In the second argument of this product, the Lipschwitz condition holds with a real $K$ for enough small $\lambda $'s, so we have that the absolute value of the substraction of the last two terms is less or equal to
$$
K\left[(f_1)_S(c(t)),\frac{1}{2}\lambda ^2\sum\limits_{i=1}^{n-1}((f_2)_S\circ c)''_i(t_i)e_i\right].
$$
Applying now the limit procedure at $\lambda \rightarrow 0$ we get the required equality.

In the second part $(f_1)_T$ and $(f_2)_T$ are real-real functions, respectively so
$$
\lim\limits_{\lambda \rightarrow 0}\frac{1}{\lambda}([(f_1)_T(c(t+\lambda )),(f_2)_T(c(t+\lambda ))]-[(f_1)_T(c(t)),(f_2)_T(c(t))])=
$$
$$
=-((f_1)_T\circ c)'(t)(f_2)_T(c(t)) -(f_1)_T(c(t))((f_2)_T\circ c)'(t).
$$
Hence we have
$$
([(f_1\circ c)(t)),(f_2\circ c)(t))]^+)'=
$$
$$
=[D((f_1)_S \circ c)(t),((f_2)_S\circ c)(t))]+[(f_1)_S(c(t)),\cdot]'_{D((f_2)_S\circ{c})(t)}(((f_2)_S\circ c)(t)))-
$$
$$
-((f_1)_T\circ c)'(t)(f_2)_T(c(t))-(f_1)_T(c(t))((f_2)_T\circ c)'(t)=
$$
$$
=[D(f_1\circ c)(t),f_2(c(t))]^+ +
$$
$$
+\left([(f_1)_S(c(t)),\cdot]'_{D((f_2)_S\circ c)(t)}((f_2)_S(c(t))) - (f_1)_T(c(t))((f_2)_T\circ c)'(t)\right),
$$
and the statement is proved.
\qed

Let $F$ be a hypersurface defined by the function $f:S\longrightarrow V$. Here $f(s)=s+\mathfrak{f}(s)e_n$ denotes the point of $F$. The curve $c:\mathbb{R}\longrightarrow S$ define  a curve on $F$. We assume that $c$ is a $C^2$-curve. The following definition is very important one.

\begin{defi}
The \emph{first fundamental form} in a point $(f(c(t))$ of the hypersurface $F$  is the product
$$
\mathrm I_{f(c(t)}:=[D(f\circ c)(t),D(f\circ c)(t)]^+.
$$
\end{defi}

The variable of the first fundamental form is a tangent vector, the tangent vector of the variable curve $c$.

We can see that it is homogeneous of the second order but (in general) it has no a bilinear representation.

In fact, by the definition of $f$, if $\{e_i : i=1\cdots n-1\}$ is a basis in $S$ then the computation
$$
\mathrm I_{f(c(t))}=[\dot{c}(t)+(\mathfrak{f}\circ c)'(t)e_n, \dot{c}(t)+(\mathfrak{f}\circ c)'(t)e_n]^+=
$$
$$
=[\dot{c}(t),\dot{c}(t)]-[(\mathfrak{f}\circ c)'(t)]^2=
[\dot{c}(t),\dot{c}(t)] -\sum\limits_{i,j=1}^{n-1}\dot{c}_i(t)\dot{c}_j(t) \mathfrak{f}'_{e_i}(c(t))\mathfrak{f}'_{e_j}(c(t))= $$
$$
=[\dot{c}(t),\dot{c}(t)]-\dot{c}(t)^T\left[\mathfrak{f}'_{e_i}(c(t)) \mathfrak{f}'_{e_j}(c(t))\right]_{i,j=1}^{n-1}\dot{c}(t)
$$
shows that it is not a quadratic form. It would be a quadratic form if and only if the quantity
$$
[\dot{c}(t),\dot{c}(t)]-\dot{c}(t)^T\dot{c}(t)= [\dot{c}(t),\dot{c}(t)]-\sum\limits_{i=1}^{n-1}\dot{c}^2_i(t)
$$
vanishes. Thus if the Minkowski product is an i.p. than we can assume that the basis $\{e_1,\ldots ,e_{n-1}\}$ in $S$ is orthonormal, the mentioned difference vanishes, and $c_i(t)=\langle e_i,c(t)\rangle=\langle c(t), e_i\rangle$ and $\dot{c}(t)=\sum \limits_{i=1}^{n-1}\dot{c}_i(t)e_i$. So
$$
\mathrm I_{f(c(t))}=\dot{c}(t)^T\left(\mathrm{Id}-\left[\mathfrak{f}'_{e_i}(c(t)) \mathfrak{f}'_{e_j}(c(t))\right]_{i,j=1}^{n-1}\right)\dot{c}(t),
$$
and we get back the classical local quadratic representation of the first fundamental form. Now if $c_i(t)=0$ for $i\geq 3$ then
$$
\det I = 1 -(\mathfrak{f}'_{e_1}(c(t)))^2 -(\mathfrak{f}'_{e_2}(c(t)))^2.
$$

We now extend the definition of the second fundamental form take into consideration that the product has neither symmetry nor bilinearity properties. If $v$ is a tangent vector and $n$ is a normal vector of the hypersurface at its point $f(c(t))$ then we have
$$
0=[v,n]^+=[D(f\circ c)(t),(f\circ c)(t)]^+.
$$
Using Lemma 4 and the notation follows it, we get
$$
0=([D(f\circ c)(t),(n\circ c)(t)]^+)'=
$$
$$
=[D^2(f\circ c),n(c(t))]^++{[D(f\circ c)(t),\cdot]^+}'_{D(n\circ c)(t)}(n(c(t))).
$$
We introduce the unit normal vector fields $n^0$ by the definition $$
n^0(c(t)):=\left\{\begin{array}{cc}n(c(t)) & \mbox{ if } n \mbox{ light-like vector}\\
\frac{n(c(t))} {\sqrt{|[n(c(t)),n(c(t))]^+|}} & \mbox{ otherwise. }\end{array}\right.
$$
\begin{defi}
The \emph{second fundamental form} at the point $f(c(t))$ defined by one of the equivalent formulas:
$$
\mathrm I \mathrm I:=[D^2(f\circ c)(t),(n^0\circ c)(t)]^+_{(f\circ c)(t)}=-{[D(f\circ c)(t),\cdot]^+}'_{D(n^0\circ c)(t)}((n^0\circ c)(t)).
$$
\end{defi}
By the structure of the generalized space-time model assuming that $n(s)=s+\mathfrak{n}(s)e_n$ we get that
$$
\mathrm I \mathrm I=[D^2(f\circ c)(t),(n^0\circ c)(t)]^+_{(f\circ c)(t)}=
$$
$$
=\left[D(\dot{c}(t)+D(\mathfrak{f}\circ c)(t)e_n),\frac{ {c}(t)+(\mathfrak{n}\circ c)(t)e_n}{\sqrt{|[c(t),c(t)]-(\mathfrak{n}(c(t)))^2|}}\right]^+=
$$
$$
=\frac{\left[\ddot{c}(t)+ \left(\dot{c}(t)^T\left[\mathfrak{f}''_{e_i,e_j}|_{c(t)}\right]\dot{c}(t) +\left[\mathfrak{f}'_{e_i}|_{c(t)}\right]\ddot{c}(t)\right)e_n,c(t)+ \mathfrak{n}(c(t))e_n\right]^+} {\sqrt{|[c(t),c(t)]-(\mathfrak{n}(c(t)))^2|}}=
$$
$$
=\frac{\left[\ddot{c}(t)+ [\mathfrak{f}'_{e_i}|_{c(t)}]\ddot{c}(t)e_n,(n\circ c)(t)\right]^+-\left(\dot{c}(t)^T\left[\mathfrak{f}''_{e_i,e_j}|_{c(t)}\right] \dot{c}(t)\right)(\mathfrak{n}(c(t))} {\sqrt{|[c(t),c(t)]-(\mathfrak{n}(c(t)))^2|}}=
$$
$$
=\frac{\left[D(f)|_{c(t)}\ddot{c}(t),(n\circ c)(t)\right]^+-\left(\dot{c}(t)^T\left[\mathfrak{f}''_{e_i,e_j}|_{c(t)}\right] \dot{c}(t)\right)(\mathfrak{n}(c(t))} {\sqrt{|[c(t),c(t)]-(\mathfrak{n}(c(t)))^2|}}=
$$
$$
=-\left(\dot{c}(t)^T\left[\frac{\mathfrak{f}''_{e_i,e_j}|_{c(t)} \mathfrak{n}(c(t))} {\sqrt{|[c(t),c(t)]-(\mathfrak{n}(c(t)))^2|}}\right]_{i,j=1}^{n-1} \dot{c}(t)\right).
$$
We now can adopt a determinant of this fundamental form. It is the determinant of its quadratic form:
$$
\det \mathrm I \mathrm I:=\det\left(\left[\frac{\mathfrak{f}''_{e_i,e_j}|_{c(t)} \mathfrak{n}(c(t))} {\sqrt{|[c(t),c(t)]-(\mathfrak{n}(c(t)))^2|}}\right]_{i,j=1}^{n-1}\right).
$$

If we consider a two-plane in the tangent hyperplane then it has a two dimensional pre-image in $S$ by the regular linear mapping $Df$. The getting plane is a normed one and we can consider an Auerbach basis $\{e_1,e_2\}$ in it.
\begin{defi}
The \emph{sectional principal curvature} of a 2-section of the tangent hyperplane in the direction of the 2-plane spanned by $\{u=Df(e_1)$ and $v=Df(e_2)\}$ are the extremal values of the function
$$
\rho (D(f\circ c)):=\frac{\mathrm I \mathrm I_{f\circ c(t)}}{\mathrm I_{f\circ c(t)}},
$$
of the variable $D(f\circ c)$. We denote them by $\rho(u,v)_{\mathrm{max}}$ and $\rho(u,v)_{\mathrm{min}}$, respectively. The \emph{sectional (Gauss) curvature} $\kappa(u,v)$  (at the examined point $c(t)$) is the product
$$
\kappa(u,v):=[n^0(c(t)),n^0(c(t))]^+ \rho (u,v)_{\mathrm{max}}\rho(u,v)_{\mathrm{min}}.
$$
\end{defi}

In the case of a symmetric and bilinear product, both of the fundamental forms are quadratic and the sectional principal curvatures attained in orthogonal directions. They are the eigenvalues of the pair of quadratic forms $\mathrm I \mathrm I_{f\circ c(t)}$ and $\mathrm I_{f\circ c(t)}$. This implies that $\rho(u,v)_{\mathrm{max}}$ and $\rho(u,v)_{\mathrm{min}}$ are the solutions of the equality:
$$
0=\det\left(\mathrm I \mathrm I_{f\circ c(t)}-\lambda \mathrm I_{f\circ c(t)}\right)=\det\left(\mathrm I_{f\circ c(t)}\right)\det\left((\mathrm I_{f\circ c(t)})^{-1}\mathrm I \mathrm I_{f\circ c(t)}-\lambda \mathrm{Id}\right),
$$
showing that
$$
\kappa(u,v):=[n^0(c(t)),n^0(c(t))]^+ \rho (u,v)_{\mathrm{max}}\rho(u,v)_{\mathrm{min}}=
$$
$$
=[n^0(c(t)),n^0(c(t))]^+ \det\left(\mathrm I_{f\circ c(t)}^{-1}\mathrm I \mathrm I_{f\circ c(t)}\right)=[n(c(t)),n(c(t))]^+ \frac{\det\mathrm I \mathrm I_{f\circ c(t)}}{\det \mathrm I_{f\circ c(t)}}=
$$
$$
=[n^0(c(t)),n^0(c(t))]^+ \frac{\left(\mathfrak{f}''_{e_1,e_1}|_{c(t)} \mathfrak{f}''_{e_2,e_2}|_{c(t)}- \left(\mathfrak{f}''_{e_1,e_2}|_{c(t)}\right)^2\right) \left(\mathfrak{n}(c(t))\right)^2}{\left(1 -(\mathfrak{f}'_{e_1}(c(t)))^2 -(\mathfrak{f}'_{e_2}(c(t)))^2\right){|[c(t),c(t)] -(\mathfrak{n}(c(t)))^2|}}.
$$
But we can choose for the function  $n$
$$
n(c(t)):=\mathfrak{f}'_{e_1}(c(t))e_1+\mathfrak{f}'_{e_2}(c(t))e_2+e_n
$$
with $\mathfrak{n}(c(t))=1$ and for a 2-plane of the tangent hyperplane which contains only space-like vectors and has time-like normal vector with absolute value
$$
[n(c(t)),n(c(t))]^+=
\sqrt{1-(\mathfrak{f}'_{e_1}(c(t)))^2 -(\mathfrak{f}'_{e_2}(c(t)))^2}
$$
getting the well-known formula
$$
\kappa(u,v)=\frac{-\mathfrak{f}''_{e_1,e_1}|_{c(t)} \mathfrak{f}''_{e_2,e_2}|_{c(t)} +\left(\mathfrak{f}''_{e_1,e_2}|_{c(t)}\right)^2 }{\left(1 -(\mathfrak{f}'_{e_1}(c(t)))^2 -(\mathfrak{f}'_{e_2}(c(t)))^2\right)^2}
$$
(see in \cite{dubrovin} p.95.).

The Ricci curvature of a Riemannian hypersurface at a point $p=(f\circ c)(t)$ in the direction of the tangent vector $v=D(f\circ c)$ is the sum of the sectional curvatures in the directions of the planes spanned by the tangent vectors $v$ and $u_i$, where $u_i$ are the vectors of an orthonormal basis of the orthogonal complement of $v$. This value is independent from the choosing of the basis. Choose random (by uniform distribution) the orthonormal basis! (\cite{balazs}) The corresponding sectional curvatures $\kappa(u_i,v)$ will be random variables with the same expected values. The sum of them is again a random variable which expected value corresponding to the Ricci curvature at $p$ with respect to $v$. Hence it is equal to $n-2$-times the expected value of the random sectional curvature determined by all of the two planes through $v$. Similarly the scalar curvature of the hypersurface at a point is the sum of the sectional curvatures defined by any two vectors of an orthonormal basis of the tangent space, it is also can be considered as an expected value. This motivates the following definition:

\begin{defi}
The \emph{Ricci curvature} $\mathrm{Ric}(v)$ in the direction $v$ at the point $f(c(t))$ is
$$
\mathrm{Ric}(v)_{f(c(t))} := (n - 2)\cdot E(\kappa_{f(c(t))}(u,v))
$$
where $\kappa_{f(c(t))}(u,v)$ is the random variable of the sectional curvatures of the two planes spanned by $v$ and a random  $u$ of the tangent hyperplane holding the equality $[u,v]^+=0$. We also say that the scalar curvature of the hypersurface $f$ at its point $f(c(t))$ is
$$
\Gamma_{f(c(t))}:= {n-1 \choose 2} \cdot E(\kappa_{f(c(t))}(u,v)).
$$
\end{defi}

\subsection{Arc-length}

In this section we also assume that the s.i.p. of $S$ is continuously differentiable.
If the first fundamental form is positive we can adopt for the curves  well-defined arc-lengthes and we can define a metric on the hypersurface.

The following definition was used in \cite{gho} for the metric of the imaginary unit sphere. We now adopt it for an arbitrary hypersurface.

\begin{defi} Denote by $p,q$ a pair of points in $F$ where $F$ is a hypersurface of the generalized space-time model. Consider the set $\Gamma _{p,q}$ of equally oriented piecewise differentiable curves $(f\circ c)(t)$,
$a\leq t\leq b$, of $F$ emanating from $p$ and terminating at $q$. Then the \emph{ pre-distance} of these points is
$$
\rho (p,q)=\inf \left\{ \int \limits _{a}^{b}\sqrt{|\mathrm I_{(f\circ c)(x)}|}dx \mbox{ for } f\circ c\in \Gamma _{p,q}\right \}.
$$
\end{defi}

It is easy to see that pre-distance satisfies the
triangle inequality; thus it gives a metric on $F$ (see \cite{tamassy}).  On a hypersurface which contains only space-like tangent vectors it is the usual definition of the Minkowski-Finsler distance. Such hypersurfaces were called by space-like ones, we mention for an example the imaginary unit sphere.
Introduce the {\bf arc-length function} $l_a(\tau)$ of a curve $f\circ c$ for which the light-like points gives a closed, zero measured set by the function
$$
l_a(\tau)=\int \limits _{a}^{\tau}\sqrt{|[D({f\circ c})(x),D(f\circ c)(x)]^+_{f(c(x))}|}dx=\int \limits _{a}^{\tau}\sqrt{|\mathrm I_{f(c(x))}|}dx.
$$
Give parameters only those points of the curve in which the tangent vector of the curve is a non-light-like one. Thus a corresponding reparametrization could be well-defined and we get a pair of inverse formulas which are almost all valid; we have that
$$
\left(l_a(\tau)\right)'=\sqrt{\mathrm |I_{f(c(\tau))}|},
$$
and for the inverse function $\tau(l_a):[0,\varepsilon)\longrightarrow [a,l_a^{-1}(\varepsilon))$ holds
$$
\left(\tau (l_a)\right)'=\left(l_a^{-1}(\tau)\right)'=\frac{1}{\sqrt{|\mathrm I_{f(c(\tau(l_a)))}}|}.
$$
\begin{theorem}
Consider a curve lying on the hypersurface determining a vector field in $V$. Using the arc-length as parameter, the absolute values of the first derivative (tangent) vectors are equal to 1, moreover the second derivative vector fields are orthogonal to the first one. (With respect to the Minkowski product of $V$.)
\end{theorem}

\proof
By definition the tangent vectors are non-light-like. Thus the required differential is
$$
D\left((f\circ c)\circ (\tau (l_a))\right)=D(f\circ c)\circ(\tau (l_a))\cdot \left(\tau (l_a)\right)'=
$$
$$
=D(f\circ c)\circ (\tau (l_a))\frac{1}{\sqrt{|\mathrm I_{f(c(\tau(l_a)))}|}},
$$
implying that
$$
[D\left((f\circ c)\circ (\tau (l_a))\right),D\left((f\circ c)\circ (\tau (l_a))\right)]^+=\frac{\mathrm I_{f(c(\tau(l_a)))}}{|\mathrm I_{f(c(\tau(l_a)))}|}=\mbox{sign}(\mathrm I_{f(c(\tau(l_a)))}).
$$
If we would like to consider the derivative of the tangent vector field, we have to compute
$$
D\left(D\left((f\circ c)\circ (\tau (l_a))\right)\right).
$$
By Lemma 4 we get
$$
D\left(\sqrt{|\mathrm I_{f(c(\tau(l_a)))}|}\right)
=
$$
$$
=\frac{\mbox{sign}(\mathrm I_{f(c(\tau(l_a)))})}{2\sqrt{|\mathrm I_{f(c(\tau(l_a)))}}|}
\left\{\frac{[D^2({f\circ c})\circ(\tau (l_a)),D(f\circ c)\circ (\tau (l_a))]}{\sqrt{|\mathrm I_{f(c(\tau(l_a)))}}|}+\right.
$$
$$
\left.+[D({f\circ c})\circ (\tau (l_a)),\cdot ]'_{D((D{f\circ c})\circ \tau (l_a))}(D(f\circ c)\circ (\tau(l_a)))\right\},
$$
and since the s.i.p is continuously differentiable we have by Theorem 5 that
$$
[D({f\circ c})\circ (\tau (l_a)),\cdot ]'_{D((D{f\circ c})\circ \tau (l_a))}(D(f\circ c)\circ (\tau(l_a)))=
$$
$$
=\frac{[D^2({f\circ c})\circ(\tau (l_a)),D(f\circ c)\circ (\tau (l_a))]}{\sqrt{|\mathrm I_{f(c(\tau(l_a)))}}|}.
$$
Thus the complete differential is
$$
D\left(D(f\circ c)\circ(\tau (l_a))\frac{1}{\sqrt{|\mathrm I_{f(c(\tau(l_a)))}|}}\right)
=\frac{D^2({f\circ c})\circ (\tau (l_a))}{|\mathrm I_{f(c(\tau(l_a)))}|}-
$$
$$
-\mbox{sign}(\mathrm I_{f(c(\tau(l_a)))})\frac{[D^2({f\circ c})\circ(\tau (l_a)),D(f\circ c)\circ(\tau (l_a))]}{(\mathrm I_{f(c(\tau(l_a)))})^2}D(f\circ c)\circ(\tau (l_a)).
$$
Now we can see that
$$
\left[D\left(D(f\circ c)\circ(\tau (l_a))\frac{1}{\sqrt{\mathrm I_{f(c(\tau(l_a)))}}}\right),D(f\circ c)\circ (\tau (l_a))\right]=0
$$
as we stated.
\qed

Inner metric determines the geodesics of the hypersurface in a standard way, as the curves representing the infinum in Definition 13. Since we use the concept of arc-length parametrization for a general (not necessary space-like) curve, we can introduce its {\bf velocity and acceleration vectors fields} as the first and second derivative vector fields of the natural parametrization, respectively. We can introduce the concept of geodesics as the solutions of the Euler-Lagrange equation with respect to the hypersurface. More precisely:
\begin{defi}
We say that the $C^2$-curve $f\circ c$ (with almost all non-light-like tangent vectors) is a \emph{geodesic} of the hypersurface $F$, if its acceleration vector field is orthogonal to the tangent hyperplane of $F$ at each point of the curve. So there exists a function $\alpha (\tau (l_a)):\mathbb{R}\longrightarrow \mathbb{R}$ such that
$$
D\left(D\left((f\circ c)\circ (\tau (l_a))\right)\right)=\alpha (\tau (l_a))(n^0\circ c)(\tau (l_a)).
$$
\end{defi}

The curvature of a curve can be defined as the square root of the absolute value of the derivative of its tangent vectors with respect to this parametrization.
\begin{defi}
The \emph{curvature of the curve} $f\circ c$ is the non-negative function
$$
\gamma_{f\circ c}(\tau(l_a)):=\sqrt{\left|[D\left(D\left((f\circ c)\circ (\tau (l_a))\right)\right),D\left(D\left((f\circ c)\circ (\tau (l_a))\right)\right)]^+\right|}=
$$
$$
=|\alpha (\tau (l_a))|.
$$
\end{defi}

If the curvature is non-zero then we can define the vector $(m\circ c)(\tau(l_a))$ by the equality:
$$
(m\circ c)(\tau(l_a))=\frac{D\left(D\left((f\circ c)\circ (\tau (l_a))\right)\right)}{\gamma_{f\circ c}(\tau(l_a))}.
$$
From this equality immediately follows that
$$
[(m\circ c)(\tau(l_a)),(m\circ c)(\tau(l_a))]^+=
$$
$$
=\frac{[D\left(D\left((f\circ c)\circ (\tau (l_a))\right)\right),D\left(D\left((f\circ c)\circ (\tau (l_a))\right)\right)]^+}{\gamma_{f\circ c}^2(\tau(l_a))}.
$$
Using the equality
$$
D\left(D\left((f\circ c)\circ (\tau (l_a))\right)\right)=D\left(D(f\circ c)\circ(\tau (l_a))\frac{1}{\sqrt{|\mathrm I_{f(c(\tau(l_a)))}|}}\right)
=
$$
$$
=\frac{D^2({f\circ c})\circ (\tau (l_a))}{|\mathrm I_{f(c(\tau(l_a)))}|}-
$$
$$
-\mbox{sign}(\mathrm I_{f(c(\tau(l_a)))})\frac{[D^2({f\circ c})\circ(\tau (l_a)),D(f\circ c)\circ(\tau (l_a))]}{(\mathrm I_{f(c(\tau(l_a)))})^2}D(f\circ c)\circ(\tau (l_a)),
$$
computed in Theorem 7, and the orthogonality property of the vectors $D(f\circ c)$ and $n^0\circ c$, we get
a connection analogous to the Meusnier's theorem:
$$
\gamma_{f\circ c}(\tau(l_a))[(m\circ c)(\tau(l_a)),(n^0\circ c)(\tau(l_a)]^+=\left[\frac{D^2({f\circ c})\circ(\tau (l_a))}{|\mathrm I_{f(c(\tau(l_a)))}|},(n^0\circ c)(\tau(l_a)\right]^+
$$
meaning that
$$
\gamma_{f\circ c}(\tau(l_a))[(m\circ c)(\tau(l_a)),(n^0\circ c)(\tau(l_a)]^+=\frac{\mathrm I \mathrm I_{f(c(\tau(l_a)))}}{|\mathrm I_{f(c(\tau(l_a)))}|}.
$$
The product form of this equality is
$$
\gamma_{f\circ c}(\tau(l_a))[(m\circ c)(\tau(l_a)),(n^0\circ c)(\tau(l_a)]^+|\mathrm I_{f(c(\tau(l_a)))}|=\mathrm I \mathrm I_{f(c(\tau(l_a)))}.
$$
This for light-like vectors is also valid, if we define their acceleration vectors as  vectors of zero length. By definition, for a geodesic curve
$$
[(m\circ c)(\tau(l_a)),(m\circ c)(\tau(l_a))]^+=\frac{(\alpha (\tau (l_a)))^2}{(\gamma_{f\circ c}(\tau(l_a))�^2}\left[(n^0\circ c)(\tau (l_a)),(n^0\circ c)(\tau (l_a))\right]^+,
$$
showing that $m\circ c$ and $n^0\circ c$ have the same casual characters and thus
$$
m\circ c=\mbox{sign}(\alpha (\tau (l_a))) (n^0\circ c).
$$
Thus the product form of the Meusnier's theorem simplified into the equality
$$
\alpha (\tau (l_a))[(n^0\circ c)(\tau(l_a)),(n^0\circ c)(\tau(l_a)]^+|\mathrm I_{f(c(\tau(l_a)))}|=\mathrm I \mathrm I_{f(c(\tau(l_a)))}.
$$
Equivalently we get
$$
\alpha (\tau (l_a))=[(n^0\circ c)(\tau(l_a)),(n^0\circ c)(\tau(l_a)]^+\frac{\mathrm I \mathrm I_{f(c(\tau(l_a)))}}{|\mathrm I_{f(c(\tau(l_a)))}|}=
$$
$$
=[(n^0\circ c)(\tau(l_a)),(n^0\circ c)(\tau(l_a)]^+\mbox{sign}(\mathrm I_{f(c(\tau(l_a)))})\rho (D(f\circ c)).
$$

If all tangent vectors are space-like vectors and the normal ones are time-like vectors, respectively, then the extremal values of the function $\alpha (\tau (l_a))=-\rho (D(f\circ c))$ on a two plane are the negatives of the principal curvatures. By the homogeneity properties of the fundamental forms, the investigated functions can be restricted to such a special subset, on which all of the possible values attain, to the unit circle of this plane. This set is compact and thus there are two extremal values  and at least two corresponding unit vectors, respectively. The convexity of such a hypersurface implies that the signs of the extremal values of $\alpha (\tau (l_a))$ are equals, so the two principal curvatures has the same signs and thus the sectional curvature is negative.

On the other hand, the characters of such a tangent plane would be only two types; either it is a space-like plane containing only space-like vectors or it has two non-paralel light-like vectors partitioning the plane two double cones one of them contains the space-like vectors and the other one the time-like vectors, respectively.

In the second case, we can restrict our function onto the union of the imaginary unit circle, the de Sitter circle, and the two lines containing the light-like vectors, respectively. We omit the two direction of the light-like vectors and we can determine the extremal values of the second fundamental form on the de Sitter sphere and on the imaginary unit sphere, respectively.
For example if the signs of the functions $\alpha (\tau (l_a))$ and $\mathrm I_{f(c(\tau(l_a)))}$ are equals, and the normal vectors are space-time vectors, then the principal curvatures have the same signs, implying that their product is positive. In this case, the sectional curvature is positive.

\section{Four interesting premanifolds}

In this section we give the most important hypersurfaces of a generalized space-time model and determine their geometries, respectively.

\subsection{Imaginary unit sphere}

Then by Theorem 6 $(H^+,ds^2)$ is a Minkowski-Finsler space, where for the vectors $u_1$ and $u_2$ of $T_v$ we have
$$
ds^2_v(u_1,u_2)=[u_1,u_2]^+_v
$$
with the Minkowski product $[\cdot,\cdot]^+_v$ of the tangent space $T_v$. This gives a possibility to examine the geometric property of $H^+$ on the base of the standard differential geometry of a space-time hypersurface.
First we prove the following theorem:
\begin{theorem}
$H^+$ is always convex. It is strictly convex if and only if the s.i.p. space $S$ is a strictly convex space.
\end{theorem}

\proof
Let $w=s'+t'$ be a point of $H^+$ and consider the product
$$
[w-v,v]^+=[s'-s,s]+[t'-t,t]=[s',s]-[s,s]-(\lambda '-\lambda)\lambda =[s',s]-\lambda' \lambda +1,
$$
where $t'=\lambda 'e_n$, $t=\lambda e_n$ and $s',s\in S$ with positive $\lambda '$ and $\lambda$, respectively.
Since
$$
\sqrt{1+[s',s']}=\lambda' \mbox{ and } \sqrt{1+[s,s]}=\lambda
$$
thus
$$
[w-v,v]^+=[s',s]-\sqrt{1+[s',s']}\sqrt{1+[s,s]}+1\leq
$$
$$
\leq \sqrt{[s',s'][s,s]}-\sqrt{1+[s',s']}\sqrt{1+[s,s]}+1\leq 0,
$$
because of the relation
$$
[s',s'][s,s]+2\sqrt{[s',s'][s,s]}+1\leq [s',s'][s,s]+([s',s']+[s,s])+1.
$$
(We used here the inequality between the arithmetic and geometric means of two positive numbers.)
Remark that equality holds if and only if the norms of $s'$ and $s$ are equal to each other and thus $\lambda'=\lambda$, too. So we have $$
[s',s]-[s,s]=0,
$$
or equivalently
$$
[s',s]=\sqrt{[s',s'][s,s]}.
$$
From the characterization of the strict convexity of an s.i.p. space we get $H^+$ contains only the point $v$ of the tangent space $T_v$ if and only if the s.i.p. space $S$ is strictly convex. \qed

To determine the first fundamental form consider the map $h=s+\mathfrak{h}(s)e_n$ giving the points of $H^+$. (Here $\mathfrak{h}(s)=\sqrt{1+[s,s]}$ is a real valued function.) Then we get that
$$
\mathrm I=[\dot{c}(t)+(\mathfrak{h}\circ c)'(t)e_n, \dot{c}(t)+(\mathfrak{h}\circ c)'(t)e_n]^+=
$$
$$
=[\dot{c}(t),\dot{c}(t)]-[(\mathfrak{h}\circ c)'(t)]^2,
$$
where $\dot{c}(t)$ means the tangent vector of the curve $c$ of $S$ at its point $c(t)$. Using Lemma 3 and Theorem 5 we have
$$
\mathrm I =[\dot{c},\dot{c}]-\frac{\left([\dot{c}(t),c(t)]+ [c(t),\cdot]'_{\dot{c}(t)}(c(t))\right)^2}{4(1+[c(t),c(t)])}= [\dot{c},\dot{c}]-\frac{[\dot{c}(t),c(t)]^2}{1+[c(t),c(t)]}.
$$

From this formula, by the Cauchy-Schwartz inequality, we can get a new proof for the fact that this form is positive.

The second fundamental form of $H^+$ is
$$
\mathrm I \mathrm I:=[\ddot{c}(t)+(\mathfrak{h}\circ c)''(t)e_n,c(t)+(\mathfrak{h}\circ c)(t)e_n]^+_{(\mathfrak{h}\circ c)(t)}=[\ddot{c}(t),c(t)]-(\mathfrak{h}\circ c)''(t)\mathfrak{h}(c(t)),
$$
since
$$
n\circ c=h\circ c=c(t)+(\mathfrak{h}\circ c)(t)e_n.
$$
First we compute the derivative of
$$
(\mathfrak{h}\circ c)'(t):\mathbb{R}\longrightarrow \mathbb{R}
$$
at its point $t$. We use again the formulas of Lemma 3 and Lemma 4 getting
$$
(\mathfrak{h}\circ c)''(t)=\left((\mathfrak{h}\circ c)'\right)'(t)=\left(\frac{[\dot{c}(t),c(t)]} {\sqrt{1+[c(t),c(t)]}}\right)'=
$$
$$
=\frac{[\dot{c}(t),c(t)]'}{\sqrt{1+[c(t),c(t)]}}-
\frac{\frac{[\dot{c}(t),c(t)]}{\sqrt{1+[c(t),c(t)]}} [\dot{c}(t),c(t)]}{(1+[c(t),c(t)])}
$$
and so
$$
(\mathfrak{h}\circ c)''(t)\mathfrak{h}(c(t))=[\dot{c}(t),c(t)]'-\frac{[\dot{c}(t),c(t)]^2}{1+[c(t),c(t)]}=
$$
$$
\left([\ddot{c}(t),c(t)]+[\dot{c}(t),\cdot ]'_{\dot{c}(t)}(c(t))\right)-\frac{[\dot{c}(t),c(t)]^2}{1+[c(t),c(t)]}.
$$
Thus the second fundamental form is
$$
\mathrm I \mathrm I=-[\dot{c}(t),\cdot ]'_{\dot{c}(t)}(c(t))+\frac{[\dot{c}(t),c(t)]^2}{1+[c(t),c(t)]},
$$
or using the formula
$$
\|y\|\|\cdot\|''_{x,z}(y)= [x,\cdot]'_z(y)-\frac{\Re[x,y]\Re[z,y]}{\|y\|^2},
$$
we have equivalently
$$
\mathrm I \mathrm I=-\|{c}(t)\|\|\cdot\|''_{\dot{c}(t),\dot{c}(t)}{c}(t)- \frac{[\dot{c}(t),c(t)]^2}{\|c(t)\|^2(1+\|c(t)\|^2)}.
$$
If we also assume that the norm is a $C^2$ function of its argument then we can use Theorem 5 and we get
$$
\mathrm I \mathrm I= -[\dot{c}(t),\dot{c}(t)]+\frac{[\dot{c}(t),c(t)]^2}{1+[c(t),c(t)]}= -\mathrm I.
$$

By the positivity of the first fundamental form on $H^+$, we get that the second fundamental form is negative definite and
$$
\rho (u,v)_{\mathrm{max}}=\rho(u,v)_{\mathrm{min}}=-1.
$$
This implies that the sectional curvatures are equal to $-1$, the Ricci and scalar curvatures in any direction at any point is $-(n-2)$ and $-{n-1 \choose 2}$, respectively. We proved:

\begin{theorem}
If the $S$ is a continuously differentiable s.i.p. space then the imaginary unit sphere has constant negative curvature.
\end{theorem}

Observe that our definitions in the case when the Minkowski product is an i.i.p. go to the usual concepts of hypersurfaces of a semi-Riemann manifolds (see \cite{moussong}, \cite{oneill} or \cite{verpoort}) so we can regard $H^+$ a natural generalization of the usual hyperbolic space. Thus we can say that $H$ is premanifold with constant negative curvature and $H^+$ is a {\bf prehyperbolic} space.

\subsection{de Sitter sphere}

In this subsection we shall investigate the hypersurface of those points of a generalized space-time model which scalar square is equal to one. In a pseudo-euclidean space this set was called by the {\bf de Sitter sphere}. The tangent hyperplanes of the de Sitter space are pseudo-euclidean spaces. We will denote by $G$ this set. $G$ is not a hypersurface of $V$ but we can restrict our investigation to the positive part of $G$ defined by
$$
G^+=\{s+t\in G \mbox{ : } t=\lambda e_n \mbox{ where } \lambda > 0\}.
$$
We remark that the local geometry of $G^+$ and $G$ is agree by the symmetry of $G$ in the subspace $S$. $G^+$ is already a hypersurface defined by the function
$$
g(s)=s+\mathfrak{g}(s)e_n,
$$
where
$$
\mathfrak{g}(s)=\sqrt{-1+[s,s]} \mbox{ for } [s,s]>1.
$$
First we calculate the directional derivatives of the function
$$
\mathfrak{g}:s\longmapsto \sqrt{-1+[s,s]}
$$
giving the corresponding tangent vectors of form
$$
u=\alpha(e+\mathfrak{g}'_e(s)e_n).
$$
Since between $\mathfrak{g}$ and $\mathfrak{f}:s\longmapsto \sqrt{1+[s,s]}$, there is the connection
$$
\mathfrak{f}^2(s)+\mathfrak{g}^2(s)=2[s,s],
$$
the derivative of $\mathfrak{g}$ in the direction of the unit vector $e\in S$ (by Lemma 1 and Lemma 3) can be calculated from the equality
$$
2\mathfrak{f}(s)\mathfrak{f}'_e(s)+2\mathfrak{g}(s)\mathfrak{g}'_e(s)=4\|s\|\|\cdot\|'_e(s)=4[e,s].
$$
Thus
$$
\mathfrak{g}'_e(s)=\frac{[e,s]}{\mathfrak{g}(s)}=\frac{[e,s]}{\sqrt{-1+[s,s]}}
$$
meaning that
$$
[u,u]^+=\alpha ^2\left(1-\frac{[e,s]^2}{(-1+[s,s])}\right)=\alpha ^2\frac{-1+[s,s]-[e,s]^2}{-1+[s,s]}.
$$
From this we can see immediately that
$$
\begin{array}{ccl}
                   &[u,u]^+    >0 & \mbox{ if } -1+[s,s]>[e,s]^2\\
                    &[u,u]^+=0 & \mbox{ if } -1+[s,s]=[e,s]^2\\
                    & [u,u]^+<0 & \mbox{ if } -1+[s,s]<[e,s]^2.
              \end{array}
$$
So a vector $s'$ of the $n-2$-subspace of $S$ orthogonal to $s$ determines a space-time tangent vector in the tangent space and a tangent vector corresponding to $\alpha s$ is a time-like one.
To determine the light-like tangent vectors consider a unit vector $e\in S$ of the form
$$
e=\frac{\pm\sqrt{-1+[s,s]}}{[s,s]} s+s', \mbox{ where } s'\in s^\bot.
$$
Such a unit vector lying in the intersection of the unit sphere of $S$ by the union of $n-2$-dimensional affine subspaces
$$
s^\bot + \frac{\pm\sqrt{-1+[s,s]}}{[s,s]} s.
$$
Since $s^\bot$ is the orthogonal complement of $s$ in $S$ and
$$
\left(\frac{\pm\sqrt{-1+[s,s]}}{[s,s]}\right)^2[s,s]= \frac{-1+[s,s]}{[s,s]}<1,
$$
this intersection is the union of two spheres of dimension $n-3$.

Thus the directions of the light-like vectors form a cone of the tangent hyperplane; the cone of the points:
$$
u=\alpha\left(\left(\pm\sqrt{-1+[s,s]}s+[s,s]s'\right)\pm [s,s]e_n\right).
$$

Recall that we considered the tangent hyperplane as a subspace of the original vector space and observe that thus we can admit it an inner Minkowskian structure, with respect to the positive and negative subspaces

$$
S':=s^{\bot}\cap S=s^{\bot} \mbox{ and } T'=\alpha \left(\sqrt{-1+[s,s]}s+[s,s]e_n\right).
$$

First we note the following:
\begin{theorem}
$G^+$ and its tangent hyperplanes are intersecting, consequently there is no point at which $G$ would be convex.
\end{theorem}
\proof
At an arbitrary point of $G^+$ there are two sets lying on $G^+$ and having in distinct halfspaces with respect to the corresponding tangent hyperplane. The first set is the intersection of the 2-plane spanned by $e_n$ and $s+t\in M$; and the other one is an arbitrary curve of the $(n-2)$-hypersurface defined by the intersection of $G$ and the hyperplane $S+(s+t)$. In fact, a normal vector of the tangent hyperplane at $s+t$ is itself $s+t$, because we have
$$
\left[e+\frac{[e,s]}{\sqrt{-1+[s,s]}}e_n,s+ {\sqrt{-1+[s,s]}}e_n\right]^+=0.
$$
Thus with $\alpha >\frac{1}{\sqrt{[s,s]}}$ we have
$$
\left[\left(\alpha s+\sqrt{-1+[\alpha s,\alpha s]}e_n\right)-\left(s+\sqrt{-1+[s,s]}e_n\right), s+\sqrt{-1+[s,s]}e_n\right]^+=
$$
$$
=(\alpha -1)[s,s]+(\sqrt{-1+[s,s]}-\sqrt{-1+[\alpha s,\alpha s]})\sqrt{-1+[s,s]}=
$$
$$
=-1+\alpha [s,s]-\sqrt{(-1+[\alpha s,\alpha s])(-1+[s,s])}=
$$
$$
=\alpha [s,s]-1-\sqrt{1-(1+\alpha ^2)[s,s]+\alpha ^2[s,s]^2}\geq 2(\alpha [s,s]-1)>2(\|s\|-1)\geq 0.
$$
On the other hand if $s'+t\in M$ arbitrary, then $\|s'\|=\|s\|$ thus $$
[s'-s+(t-t),s+t]^+=[s',s]-[s,s]\leq \sqrt{[s',s']}\sqrt{[s,s]}-[s,s]=0,
$$
with equality if and only if $s'=\pm s$.
\qed

Continue our investigation with the computation of the fundamental forms.
Using the function $g$ the first fundamental form has the form
$$
\mathrm I=[\dot{c}(t)+(\mathfrak{g}\circ c)'(t)e_n, \dot{c}(t)+(\mathfrak{g}\circ c)'(t)e_n]^+=
$$
$$
=[\dot{c}(t),\dot{c}(t)]-[(\mathfrak{g}\circ c)'(t)]^2.
$$
Using Lemma 3 and Theorem 5 we get
$$
\mathrm I =[\dot{c},\dot{c}]-\frac{\left([\dot{c}(t),c(t)]+ [c(t),\cdot]'_{\dot{c}(t)}(c(t))\right)^2}{4(-1+[c(t),c(t)])}= [\dot{c},\dot{c}]-\frac{[\dot{c}(t),c(t)]^2}{-1+[c(t),c(t)]}.
$$
Furthermore we also have that
$$
n\circ c=g\circ c=c(t)+(\mathfrak{g}\circ c)(t)e_n
$$
thus
$$
\mathrm I \mathrm I:=[\ddot{c}(t)+(\mathfrak{g}\circ c)''(t)e_n,c(t)+(\mathfrak{g}\circ c)(t)e_n]^+_{(\mathfrak{g}\circ c)(t)}=[\ddot{c}(t),c(t)]-(\mathfrak{g}\circ c)''(t)\mathfrak{g}(c(t)).
$$
The derivative of the real function
$$
(\mathfrak{g}\circ c)'(t)=D(\mathfrak{g}\circ c)(t):\mathbb{R}\longrightarrow \mathbb{R}
$$
at its point $t$ is:
$$
(\mathfrak{g}\circ c)''(t)=\frac{[\dot{c}(t),c(t)]'}{\sqrt{-1+[c(t),c(t)]}}-
\frac{\frac{[\dot{c}(t),c(t)]}{\sqrt{-1+[c(t),c(t)]}} [\dot{c}(t),c(t)]}{(-1+[c(t),c(t)])}
$$
so by Lemma 4
$$
(\mathfrak{g}\circ c)''(t)\mathfrak{g}(c(t))=[\dot{c}(t),c(t)]'- \frac{[\dot{c}(t),c(t)]^2}{-1+[c(t),c(t)]}=
$$
$$
=\left([\ddot{c}(t),c(t)]+[\dot{c}(t),\cdot ]'_{\dot{c}(t)}(c(t))\right)-\frac{[\dot{c}(t),c(t)]^2}{-1+[c(t),c(t)]}.
$$
Thus we have
$$
\mathrm I \mathrm I=-[\dot{c}(t),\cdot ]'_{\dot{c}(t)}(c(t))+\frac{[\dot{c}(t),c(t)]^2}{-1+[c(t),c(t)]}.
$$
If we assume again that the norm is a $C^2$ function of its argument then we can use again Theorem 5 and we get
$$
\mathrm I \mathrm I= -[\dot{c}(t),\dot{c}(t)]+\frac{[\dot{c}(t),c(t)]^2}{-1+[c(t),c(t)]}= -\mathrm I,
$$
as in the case of $H^+$. The principal curvatures are equal to $-1$. But the scalar squares of the normal vectors is positive at all points of $G^+$ implying that the sectional curvatures are equal to $1$. The Ricci curvatures in any directions and at any points are equal to $(n-2)$, moreover the scalar curvatures at any points are equal to ${n-1 \choose 2}$ showing that:

\begin{theorem}
The de Sitter sphere $G$ has constant positive curvature if $S$ is a continuously differentiable s.i.p space.
\end{theorem}

On the basis of this theorem we can say about $G$ as a premanifold of constant positive curvature and we may say that it is a {\bf pre-sphere}.

\subsection{The light cone}

The inner geometry of the light cone $L$ can be determined, too. Let $L^+$ be the positive part of this double cone determined by the function:
$$
l(s)=s+\sqrt{[s,s]}e_n.
$$
If $S$ is a uniformly continuous s.i.p. space, then the tangent vectors at $s$ are of the form:
$$
u=\alpha\left(e+\|\cdot\|'_e(s)e_n\right)=\alpha\left(e+ \frac{[e,s]}{\sqrt{[s,s]}}e_n\right).
$$
Thus all tangents orthogonal to $l(s)$ which is also a tangent vector. (Choose $e=s^0$ and $\alpha=\|s\|$!) But the orthogonal companion of a neutral (isotropic or light-like) vector in a s.i.i.p space is an $(n-1)$-dimensional degenerated subspace containing it (\cite{gho} (Theorem 7)) Tangent hyperplanes are exist at every points of $L^+$ and it is an $(n-1)$-dimensional degenerated subspace of $V$. This also a support hyperplane of $L$. In fact, by $v=s+t$ and $w=s'+t'$ we get
$$
[w-v,v]^+=[s',s]+[t',t]=[s',s]-\lambda '\lambda
$$
where $t'=\lambda 'e_n$, $t=\lambda e_n$ and $s',s\in S$ with positive $\lambda '$ and $\lambda$, respectively.
Since
$$
\sqrt{[s',s']}=\lambda' \mbox{ and } \sqrt{[s,s]}=\lambda
$$
thus
$$
[w-v,v]^+=[s',s]-\sqrt{[s',s']}\sqrt{[s,s]}\leq 0
$$
holds by the Cauchy-Schwartz inequality.
We remark that equality holds if and only if $s'=\alpha s$ meaning that there is only one line of $L^+$ in the tangent space $T_v$. Thus the light cone is convex and thus the second fundamental form is semi-definite quadratic form. It also follows that any other vectors of the tangent hyperplane are space-like ones and there are two types of tangent 2-planes; one of them space-like plane and the other one contains space-like vectors and a doubled line of light-like vectors. In the first case, the corresponding principal and sectional curvatures is well defined and have negative values, respectively. To determine it we compute the fundamental forms.

In the case when $S$ is continuously differentiable, the first fundamental form is
$$
\mathrm I =[\dot{c},\dot{c}]-\frac{\left([\dot{c}(t),c(t)]+ [c(t),\cdot]'_{\dot{c}(t)}(c(t))\right)^2}{4[c(t),c(t)]}= [\dot{c},\dot{c}]-\frac{[\dot{c}(t),c(t)]^2}{[c(t),c(t)]},
$$
and the second one is
$$
\mathrm I \mathrm I=-[\dot{c}(t),\cdot ]'_{\dot{c}(t)}(c(t))+\frac{[\dot{c}(t),c(t)]^2}{[c(t),c(t)]}= -[\dot{c}(t),\dot{c}(t)]+\frac{[\dot{c}(t),c(t)]^2}{[c(t),c(t)]}=-\mathrm I.
$$
Thus the principal curvatures are $-1$ as in the cases of the unit spheres. However our definition gives at such a point zero sectional curvature for it, because of the zero lengthes of the normal vectors. The above computation can be used in the second case, too. Agreed that we calculate the fundamental forms only non-light-like directions, so on the plane of the second type the principal curvatures are also $-1$ and the sectional curvatures are zero, too. This implies that the Ricci and scalar curvatures are also zero, respectively. We have got

\begin{theorem}
The light cone $L^+$ has zero curvature if $S$ is a continuously differentiable s.i.p space.
\end{theorem}

Hence $L$ is a premanifold with zero curvature and we may say that it is a {\bf pre-Euclidean} space.

\subsection{The unit sphere of the s.i.p. space $(V,[\cdot,\cdot]^-)$}

In this subsection we shall investigate the hypersurface of those points of the generalized space-time model which collects the unit sphere of the embedding s.i.p. space. In a pseudo-euclidean space it is the unit sphere of the embedding euclidean space. Its tangent hyperplanes are pseudo-euclidean one. We will denote by $K$ this set. $K$ is not a hypersurface but we can restrict our investigation to the positive part of $K$ defined by
$$
K^+=\{s+t\in K \mbox{ : } t=\lambda e_n \mbox{ where } \lambda > 0\}.
$$
$K^+$ is a hypersurface defined by the function
$$
k(s)=s+\mathfrak{k}(s)e_n,
$$
where
$$
\mathfrak{k}(s)=\sqrt{1-[s,s]} \mbox{ for } [s,s]<1.
$$
The directional derivatives of the function
$$
\mathfrak{k}:s\longmapsto \sqrt{1-[s,s]} \mbox{ for } [s,s]<1
$$
gives  the corresponding tangent vectors of form
$$
u=\alpha(e+\mathfrak{k}'_e(s)e_n).
$$
Since by the function
$$
\mathfrak{f}:s\longmapsto \sqrt{1+[s,s]},
$$
we have the equality
$$
\mathfrak{f}^2(s)+\mathfrak{k}^2(s)=2
$$
the derivative in the direction of the unit vector $e\in S$ is
$$
\mathfrak{k}'_e(s)=-\frac{[e,s]}{\sqrt{1-[s,s]}}
$$
meaning that
$$
[u,u]^+=\alpha ^2\left(1-\frac{[e,s]^2}{(1-[s,s])}\right)=\alpha ^2\frac{1-[s,s]-[e,s]^2}{1-[s,s]}.
$$
From this we can see immediately that
$$
\begin{array}{ccl}
                   &[u,u]^+    >0 & \mbox{ if } 1-[s,s]>[e,s]^2\\
                    &[u,u]^+=0 & \mbox{ if } 1-[s,s]=[e,s]^2\\
                    & [u,u]^+<0 & \mbox{ if } 1-[s,s]<[e,s]^2.
              \end{array}
$$
It follows that the vector $s'$ of the $n-2$-subspace of $S$ orthogonal to $s$ gives a space-time tangent vector and the vector corresponding to $\alpha s$ is a time-like one.

As in the case of the imaginary unit sphere we note the following:
\begin{theorem}
$K^+$ is convex. If $S$ is a strictly convex space, then $K^+$ is also strictly convex.
\end{theorem}
\proof
Let $w=s'+t'$ be a point of $K^+$ and consider the product
$$
[w-v,n_v]^+=[s'-s,s'']+[t'-t,t'']=[s',s'']-[s,s'']-(\lambda '-\lambda)\lambda'',
$$
where $t''=\lambda ''e_n$, $t'=\lambda 'e_n$, $t=\lambda e_n$ and $s'',s',s\in S$ with positive $\lambda ''$, $\lambda '$ and $\lambda$, respectively.
Since
$$
\sqrt{1-[s',s']}=\lambda' \mbox{ and } \sqrt{1-[s,s]}=\lambda
$$
and
$$
n_v=s-\sqrt{1-[s,s]}e_n
$$
thus
$$
[w-v,n_v]^+=[s',s]+\sqrt{1-[s',s']}\sqrt{1-[s,s]}-1\leq
$$
$$
\leq \sqrt{[s',s'][s,s]}+\sqrt{1-[s',s']}\sqrt{1-[s,s]}-1\leq 0,
$$
because
$$
2\sqrt{[s',s'][s,s]}\leq [s',s']+[s,s]).
$$
We remark that equality holds in the inequalities if and only if the norms of $s'$ and $s$ are equal to each other. So we have the equality
$$
[s',s]-[s,s]=0,
$$
or equivalently
$$
[s',s]=\sqrt{[s',s'][s,s]}.
$$
We also get that $v$ is the only point of $H^+$ lying on the tangent space $T_v$ if and only if the s.i.p. space $S$ is strictly convex.
\qed

Using the function $k$ the first fundamental form has the form
$$
\mathrm I=[\dot{c}(t),\dot{c}(t)]-[(\mathfrak{k}\circ c)'(t)]^2.
$$
Using Lemma 3 and Theorem 5 we have
$$
\mathrm I =[\dot{c},\dot{c}]-\frac{\left([\dot{c}(t),c(t)]+ [c(t),\cdot]'_{\dot{c}(t)}(c(t))\right)^2}{4(1-[c(t),c(t)])}= [\dot{c},\dot{c}]-\frac{[\dot{c}(t),c(t)]^2}{1-[c(t),c(t)]},
$$
and assuming that $2[c(t),c(t)]\neq 1$ we get
$$
\mathrm I \mathrm I=\left[\ddot{c}(t)+(\mathfrak{k}\circ c)''(t)e_n,\frac{c(t)-(\mathfrak{k}\circ c)(t)e_n}{\sqrt{|-1+2[c(t),c(t)]|}}\right]^+_{(\mathfrak{k}\circ c)(t)}=
$$
$$
=\frac{1}{\sqrt{|-1+2[c(t),c(t)]|}}\left([\ddot{c}(t),c(t)]+(\mathfrak{k}\circ c)''(t)\mathfrak{k}(c(t))\right).
$$
Lemma 4 implies that
$$
(\mathfrak{k}\circ c)''(t)\mathfrak{k}(c(t))=-[\dot{c}(t),c(t)]'+ \frac{[\dot{c}(t),c(t)]^2}{1-[c(t),c(t)]}=
$$
$$
=-\left([\ddot{c}(t),c(t)]+[\dot{c}(t),\cdot ]'_{\dot{c}(t)}(c(t))\right)+\frac{[\dot{c}(t),c(t)]^2}{1-[c(t),c(t)]}.
$$
thus we have
$$
\mathrm I \mathrm I=\frac{1}{\sqrt{|-1+2[c(t),c(t)]|}}\left(-[\dot{c}(t),\cdot ]'_{\dot{c}(t)}(c(t))+\frac{[\dot{c}(t),c(t)]^2}{1-[c(t),c(t)]}\right).
$$
Assuming that $S$ is continuously differentiable and using Theorem 5 we get
$$
\mathrm I \mathrm I=\frac{1}{\sqrt{|-1+2[c(t),c(t)]|}}\left( -[\dot{c}(t),\dot{c}(t)]+\frac{[\dot{c}(t),c(t)]^2}{-1+[c(t),c(t)]}\right)=
$$
$$
= -\frac{1}{\sqrt{|-1+2[c(t),c(t)]|}}\mathrm I.
$$
The principal curvatures at a point $k(c(t))$ are
$$
\rho_{\max}(u,v)=\rho_{\min}(u,v)=-\frac{1}{\sqrt{|-1+2[c(t),c(t)]|}}
$$
giving the sectional curvatures
$$
\kappa(u,v):=[n^0(c(t)),n^0(c(t))]^+ \rho (u,v)_{\mathrm{max}}\rho(u,v)_{\mathrm{min}}=
\frac{1}{-1+2[c(t),c(t)]}.
$$

The Ricci curvatures in any directions at the point $k(c(t))$ are equal to
$$
\mathrm{Ric}(v)_{k(c(t))} := (n - 2)\cdot E(\kappa_{k(c(t))}(u,v))=\frac{n-2}{-1+2[c(t),c(t)]}
$$
and the scalar curvature of the hypersurface $K^+$ at its point $k(c(t))$ is
$$
\Gamma_{k(c(t))}:= {n-1 \choose 2} \cdot E(\kappa_{f(c(t))}(u,v))=
\frac{{n-1 \choose 2}}{-1+2[c(t),c(t)]}.
$$
Finally we remark that at the points of $K^+$ having the equality
$2[c(t),c(t)]=1$ all of the curvatures can be defined as in the case of the light cone and can be regard to zero.

\begin{center}
\'Akos G. Horv\'ath,\\
 Department of Geometry \\
 Mathematical Institute \\
Budapest University of Technology and Economics\\
1521 Budapest, Hungary
\\e-mail: ghorvath@math.bme.hu
\end{center}

\end{document}